\documentclass[10pt,twoside,reqno]{amsart}
\usepackage{amssymb}
\calclayout

\numberwithin{equation}{section}
\makeatletter

\renewcommand{\@secnumfont}{\bfseries}

\renewcommand{\section}{\@startsection{section}{1}%
  {0mm}{.7\linespacing\@plus\linespacing}{.5\linespacing}
  {\normalfont\bfseries\centering}}

\newcommand{\bibsection}{\@startsection{section}{1}%
  {0mm}{.7\linespacing\@plus\linespacing}{.5\linespacing}
  {\normalfont\scshape\centering}}

\renewcommand{\@biblabel}[1]{#1.}

\newtheorem{thm}{\bf Theorem}[section]

\begin{document}

\vspace{1.3cm}

\title
        {Differential equations for Changhee  polynomials and their applications}

\author{Dmitry V. Dolgy, Dae San Kim, Taekyun Kim, Jong Jin Seo}

\thanks{\scriptsize }
\address{1\\Institute of Mathematics and Computer Science \\
           Far Eastern Federal University \\
          690950 Vladivostok, Russia}
\email{dvdolgy@gmail.com}
\address{2\\Department of  Mathematics\\
           Sogang University\\
           Seoul 121-742, Republic of Korea}
\email{dskim@sogang.ac.kr}
\address{3\\Department of  Mathematics\\
            Kwangwoon  University\\
           Seoul 139-701, Republic of Korea}
\email{tkkim@kw.ac.kr}
\address{4\\Department of Applied Mathematics\\
            Pukyong National University\\
            Busan 608-737, Republic of Korea}
\email{seo2011@pknu.ac.kr}


\maketitle

\begin{abstract} Recently, the non-linear Changhee differential equations were introduced in [5] and these differential equations turned out to be very useful for studying special polynomials and mathematical physics. Some interesting identities and properties of Changhee polynomials can also be derived from umbral calculus (\textnormal{see}\ [7]). In this paper, we consider differential equations arising from Changhee polynomials and derive some new and explicit formulae and identities from our differential equations.
\end{abstract}

\pagestyle{myheadings}
\markboth{\centerline{\scriptsize Dmitry V. Dolgy, Dae San Kim, Taekyun Kim, Jong Jin Seo}}
          {\centerline{\scriptsize Differential equations for Changhee  polynomials and their applications}}

\bigskip
\medskip
\medskip
\section{\bf Introduction}
\medskip
 As is well known, the Euler polynomials are defined by the generating function
\begin{equation}
\frac{2}{e^t +1}e^{xt} = \sum_{n=0}^{\infty}E_n (x)\frac{t^n}{n!},\ (\textnormal{see}\ [1, 3, 10]).
\end{equation}
With the viewpoint of deformed Euler polynomials, the Changhee polynomials are defined by the generating function
\begin{equation}
\frac{2}{2 +t}(1+t)^{x} = \sum_{n=0}^{\infty}Ch_n (x)\frac{t^n}{n!},\ (\textnormal{see}\ [3]).
\end{equation}
From (1.2), we note that
\begin{equation}\begin{split}
\frac{2}{e^{log(1+t)}+1} e^{xlog(1+t)} &= \sum_{n=0}^{\infty}E_n (x)\frac{1}{n!}(log(1+t))^n\\
 & =\sum_{n=0}^{\infty}E_n(x)\sum_{m=n}^{\infty}S_1(m,n)\frac{t^m}{m!}\\
 &=\sum_{m=0}^{\infty}\left(\sum_{n=0}^{m}E_n(x)S_1(m,n)\right)\frac{t^m}{m!},\\
\end{split}\end{equation}
where $S_1(m,n)$ is the stirling number of the first kind which is defined as
\begin{equation}
(x)_0=1, (x)_n=x(x-1)\cdot\cdot\cdot (x-n+1)=\sum_{l=0}^{n}S_1(n,l)x^l,  (n\geq 1).
\end{equation}
From (1.2) and (1.3), we note that
\begin{equation}
Ch_m(x)=\sum_{n=0}^{m}E_n(x)S_1(m,n),  (m\geq 0), (\textnormal{see}\ [3, 12, 13]).
\end{equation}
Replacing $t$ by $e^t-1$ in (1.2), we get
\begin{equation}\begin{split}
\sum_{n=0}^{\infty}E_n(x)\frac{t^n}{n!}=\frac{2}{e^t+1}e^{xt}&=\sum_{m=0}^{\infty}Ch_m(x)\frac{1}{m!}(e^t-1)^m\\
 &=\sum_{m=0}^{\infty}Ch_m(x)\sum_{n=m}^{\infty}S_2(n,m)\frac{t^n}{n!}\\
 &=\sum_{n=0}^{\infty}\left(\sum_{m=0}^{n}  Ch_m(x)S_2(n,m)\right)\frac{t^n}{n!},
\end{split}\end{equation}
where $S_2(n,m)$ is the stirling number of the second kind which is given by $x^n=\sum_{l=0}^{n}S_2(n,l)(x)_l, (n\geq 0).$
Thus, by (1.6), we get
\begin{equation}
E_n(x)=\sum_{m=0}^{n}Ch_m(x)S_2(n,m),  (n\geq 0), (\textnormal{see}\ [3]).
\end{equation}
Recently, several authors have studied Changhee polynomials and numbers (\textnormal{see}\ [1-20]). In this paper, we consider differential equations derived from the generating function of Changhee polynomials and give some new and explicit formulae for the Changhee polynomials by using our results on differential equations.

\section{\bf Differential equations for Chagnhee polynomials}

\medskip
Let
\begin{equation}
F=F(t,x)=\frac{1}{2+t}(1+t)^x.
\end{equation}
From (2.1), we can derive the following equations:
\begin{equation}
F^{(1)}=\frac{d}{dt}F(t,x)=(-(2+t)^{-1}+x(1+t)^{-1})F,
\end{equation}
\begin{equation}
F^{(2)}=\frac{d}{dt}F^{(1)}=(2(2+t)^{-2}-2x(2+t)^{-1}(1+t)^{-1}+(x^2-x)(1+t)^{-2})F,
\end{equation}
and
\begin{equation}\begin{split}
F^{(3)}=\frac{d}{dt}F^{(2)}&=(-6(2+t)^{-3}+6x(2+t)^{-2}(1+t)^{-1}\\
&  +(3x-3x^2)(2+t)^{-1}(1+t)^{-2}+(x^3-3x^2+2x)(1+t)^{-3})F.
\end{split}\end{equation}
Continuing this process, we put
\begin{equation}
F^{(N)}=\left(\frac{d}{dt}\right)^N F(t,x)=(\sum_{i=0}^{N}a_i(N,x)(1+t)^{-i}(2+t)^{i-N})F,
\end{equation}
where $N=0, 1, 2, \cdot \cdot \cdot .$
From (2.5), we note that
\begin{equation}\begin{split}
F^{(N+1)}=\frac{d}{dt} F^{(N)}&=\left(\sum_{i=0}^{N}a_i(N,x)(-i)(1+t)^{-i-1}(2+t)^{i-N}\right)F\\
            & \ \ \ \ \   + \left(\sum_{i=0}^{N}a_i(N,x)(1+t)^{-i}(i-N)(2+t)^{i-N-1}\right)F \\
            &\ \ \ \ \   +\left(\sum_{i=0}^{N}a_i(N,x)(1+t)^{-i}(2+t)^{i-N}\right)F^{(1)} \\
         &= \left\{\sum_{i=0}^{N}(-i) a_i(N,x)(1+t)^{-i-1}(2+t)^{i-N} \right. \\
         & \ \ \ \ \  + \sum_{i=0}^{N}(i-N)a_i(N,x)(1+t)^{-i}(2+t)^{i-N-1} \\
         & \ \ \ \ \ \  -\sum_{i=0}^{N}a_i(N,x)(1+t)^{-i}(2+t)^{i-N-1}\\
         & \ \ \ \ \ \ + \left. \sum_{i=0}^{N}xa_i(N,x)(1+t)^{-i-1}(2+t)^{i-N}   \right\}F
         \end{split}\end{equation}
         \begin{equation*}\begin{split}
         &=\left\{\sum_{i=1}^{N+1}(1-i)a_{i-1}(N,x)(1+t)^{-i}(2+t)^{i-N-1} \right.\\
         & \ \ \ \ \ \ +\sum_{i=0}^{N}(i-N)a_i(N,x)(1+t)^{-i}(2+t)^{i-N-1} \\
         & \ \ \ \ \ \ -\sum_{i=0}^{N}a_i(N,x)(1+t)^{-i}(2+t)^{i-N-1}\\
         & \ \ \ \ \ \ \left. +\sum_{i=1}^{N+1}xa_{i-1}(N,x)(1+t)^{-i}(2+t)^{i-N-1} \right\}F.
\end{split}\end{equation*}
On the other hand, by replacing $N$ by $N+1$ in (2.5), we get
\begin{equation}
F^{(N+1)}=\{\sum_{i=0}^{N+1}a_i(N+1,x)(1+t)^{-i}(2+t)^{i-N-1}\}F.
\end{equation}
Comparing the coefficients on both sides of (2.6) and (2.7), we have
\begin{equation}
a_0(N+1,x)=-(N+1)a_{0}(N,x),
\end{equation}
\begin{equation}
a_{N+1}(N+1,x)=(x-N)a_{N}(N,x),
\end{equation}
and
\begin{equation}
a_i(N+1,x)=(x+1-i)a_{i-1}(N,x)+(i-N-1)a_i(N,x),
\end{equation}
where $1\leq i \leq N.$\\
We also note that
\begin{equation}
F=F^{(0)}=a_0(0,x)F.
\end{equation}
Thus, by (2.11), we get
\begin{equation}
a_0(0,x)=1.
\end{equation}
From (2.2) and (2.5), we can derive the following equation:
\begin{equation}\begin{split}
(-(2+t)^{-1}+x(1+t)^{-1})F &=F^{(1)}=\left(\sum_{i=0}^{1}a_i(1,x)(1+t)^{-i}(2+t)^{i-1}\right)F\\
&=\left( a_0(1,x)(2+t)^{-1}+a_1(1,x)(1+t)^{-1}\right)F.
\end{split}\end{equation}
Comparing the coefficients on both sides of (2.13), we get
\begin{equation}
a_0(1,x)=-1, a_1(1,x)=x.
\end{equation}
Also, by (2.8) and (2.9), we have
\begin{equation}\begin{split}
a_0(N+1,x)&=-(N+1)a_0(N,x)=(-1)^2(N+1)Na_0(N-1,x)= \cdot \cdot\cdot\\
          &=(-1)^N(N+1)N\cdot \cdot \cdot 2a_0(1,x)=(-1)^{N+1}(N+1)!,
\end{split}\end{equation}
and
\begin{equation}\begin{split}
a_{N+1}(N+1,x)&=(x-N)a_N(N,x)=(x-N)(x-(N-1))a_{N-1}(N-1,x)\\
            &= \cdot \cdot\cdot=(x-N)(x-(N-1))\cdot \cdot \cdot (x-1)a_1(1,x)\\
          &=(x-N)(x-(N-1))\cdot \cdot\cdot(x-1)x =(x)_{N+1}.
\end{split}\end{equation}
From (2.10), we can derive the following equations:
\begin{equation}\begin{split}
a_{1}(N+1,x)&=xa_0(N,x)-Na_{1}(N,x)\\
          &=x\left(a_0(N,x)-Na_0(N-1,x) \right)+(-1)^2N(N-1)a_1(N-1,x)\\
          &=\cdot \cdot\cdot\\
          &=x\sum_{i=0}^{N-1}(-1)^i(N)_{i}a_0(N-i,x)+(-1)^NN!a_1(1,x)\\
          &=x\sum_{i=0}^{N}(-1)^i(N)_ia_0(N-i,x),
\end{split}\end{equation}
\begin{equation}\begin{split}
a_{2}(N+1,x)&=(x-1)a_1(N,x)+(1-N)a_{2}(N,x)\\
          &=(x-1)\left(a_1(N,x)+(-1)(N-1)a_1(N-1,x) \right)\\
           & \ \ \ \ \ +(-1)^2(N-1)(N-2)a_2(N-1,x)\\
          &=\cdot \cdot\cdot\\
          &=(x-1)\sum_{i=0}^{N-2}(-1)^i(N-1)_{i}a_1(N-i,x) +(-1)^{N-1}(N-1)!a_2(2,x)\\
          &=(x-1)\sum_{i=0}^{N-1}(-1)^i(N-1)_ia_1(N-i,x),
\end{split}\end{equation}
and
\begin{equation}\begin{split}
a_{3}(N+1,x)&=(x-2)a_2(N,x)+(2-N)a_{3}(N,x)\\
          &=(x-2)\left(a_2(N,x)+(-1)(N-2)a_2(N-1,x) \right)\\
          & \ \ \ \ \ +(-1)^2(N-2)(N-3)a_3(N-1,x)\\
          &=\cdot \cdot\cdot\\
\end{split}\end{equation}
\begin{equation*}\begin{split}
&=(x-2)\sum_{i=0}^{N-3}(-1)^i(N-2)_{i}a_2(N-i,x) +(-1)^{N-2}(N-2)!a_3(3,x)\\
          &=(x-2)\sum_{i=0}^{N-2}(-1)^i(N-2)_i a_2(N-i,x).
\end{split}\end{equation*}
Continuing this process, we have
\begin{equation}
a_j(N+1,x)=(x-j+1)\sum_{i=0}^{N-j+1}(-1)^i(N-j+1)_{i}a_{j-1}(N-i,x), (1 \leq j \leq N).
\end{equation}
The matrix $(a_i(j,x))_{0\leq i, j \leq N}$ is given by\\

$ \bordermatrix{
          &0        &1          &2          &\ldots  &N          \cr
    0     &1        &{-1}!      &(-1)^22!   &\ldots  &(-1)^N N!  \cr
    1     &0        &x          &\ldots     &\ldots  &\cdot      \cr
    2     &0        &0          &(x)_2      &\ldots  &\cdot      \cr
\vdots    &\vdots   &\vdots     &\vdots     &\ddots  &\vdots     \cr
    N     &0        &0          &0          &\ldots  &(x)_n      \cr}$\\

\medskip

Now, we give explicit expressions for $a_i(j,x).$ From (2.15), (2.17), (2.18), (2.19) and (2.20), we can derive the following equations:
\begin{equation}
a_1(N+1,x)=x\sum_{i=0}^{N}(-1)^i(N)_{i} a_{0}(N-i,x)=x(-1)^N(N+1)!,
\end{equation}
\begin{equation}\begin{split}
a_2(N+1,x)&=(x-1)\sum_{i_1=0}^{N-1}(-1)^{i_1}(N-1)_{i_1} a_{1}(N-i_{1},x)\\
          & =(x)_2 (-1)^{N-1}(N-1)!\sum_{i_1=0}^{N-1}(N-i_{1}),
\end{split}\end{equation}
\begin{equation}\begin{split}
a_3(N+1,x)&=(x-2)\sum_{i_2=0}^{N-2}(-1)^{i_2}(N-2)_{i_2} a_{2}(N-i_{2},x)\\
          & =(x)_3 (-1)^{N-2}(N-2)!\sum_{i_2=0}^{N-2}\sum_{i_1=0}^{N-2-i_2}(N-i_{2}-i_{1}-1),
\end{split}\end{equation}
and
\begin{equation}\begin{split}
a_4(N+1,x)&=(x-3)\sum_{i_3=0}^{N-3}(-1)^{i_3}(N-3)_{i_3} a_{3}(N-i_{3},x)\\
          & =(x)_4 (-1)^{N-3}(N-3)!\sum_{i_3=0}^{N-3}\sum_{i_2=0}^{N-3-i_3}\sum_{i_1=0}^{N-3-i_3-i_2}(N-i_{3}-i_{2}-i_{1}-2).
\end{split}\end{equation}
Continuing this process, we get
\begin{equation}\begin{split}
a_j(N+1,x)&=(x)_j(-1)^{N-j+1}(N-j+1)!\\
& \times \sum_{i_{j-1}=0}^{N-j+1}\sum_{i_{j-2}=0}^{N-j+1-i_{j-1}}
 \cdot \cdot \cdot\sum_{i_1=0}^{N-j+1-i_{j-1}-\cdot \cdot \cdot -i_2}(N-i_{j-1}\cdot \cdot\cdot  -i_{1}-j+2),
\end{split}\end{equation}
where $1 \leq j \leq N+1.$\\
Therefore, by (2.25), we obtain the following theorem.

\begin{thm}\label{Theorem 2.1}
For $N=0, 1, 2, \cdot \cdot \cdot ,$ the linear differential equations
\begin{equation*}\
F^{(N)}=\left(\sum_{i=0}^{N}a_i(N,x)(1+t)^{-i}(2+t)^{i-N} \right)F
\end{equation*}
has a solution $F=F(t,x)=\frac{1}{2+t}(1+t)^x,$\\
where
\begin{equation*}\begin{split}
&a_0(N,x)=(-1)^NN!,\\
&a_j(N,x) =(x)_j(-1)^{N-j}(N-j)! \\
 & \ \ \ \times \sum_{i_{j-1}=0}^{N-j}\sum_{i_{j-2}=0}^{N-j-i_{j-1}}
 \cdot \cdot \cdot\sum_{i_1=0}^{N-j-i_{j-1}-\cdot \cdot \cdot -i_2}(N-i_{j-1}\cdot \cdot\cdot  -i_{1}-j+1), (1\leq j \leq N).
\end{split}\end{equation*}
\end{thm}
We recall that Changhee polynomials, $Ch_n(x), (n\geq 0)$, are given by the generating function
\begin{equation}
2F=2F(t,x)=\frac{2}{2+t}(1+t)^x=\sum_{n=0}^{\infty}Ch_n(x)\frac{t^n}{n!}.
\end{equation}
On the one hand, by (2.26), we get
\begin{equation}
2F^{(N)}=\sum_{k=N}^{\infty}Ch_k(x)(k)_N\frac{t^{k-N}}{k!}=\sum_{k=0}^{\infty}Ch_{N+k}(x)\frac{t^k}{k!}.
\end{equation}
On the other hand, by Theorem 1, we have
\begin{equation}\begin{split}
2F^{(N)}&=\left( 2\sum_{i=0}^{N}a_i(N,x)(1+t)^{-i}(2+t)^{i-N}\right)F\\
        &=\sum_{i=0}^{N}2a_i(N,x)\left(\sum_{l=0}^{\infty}(-1)^l \binom{i+l-1}{l} t^l \right)\\
        & \ \ \ \times \left(\sum_{m=0}^{\infty} (-1)^m2^{i-N-m}  \binom{N+m-i-1}{m} t^m \right)\left( \sum_{p=0}^{\infty}Ch_p(x)\frac{t^p}{p!}\right)\\
        &=\sum_{i=0}^{N}a_i(N,x) \sum_{k=0}^{\infty} \sum_{l+m+p=k}^{}(-1)^{l+m}2^{i-N-m+1} \binom{i+l-1}{l} \\
        & \ \ \ \ \ \ \ \ \ \  \times \binom{N+m-i-1}{m} \frac{1}{p!}Ch_p(x)t^k\\
        &=\sum_{k=0}^{\infty}\left\{k!\sum_{i=0}^{N}a_i(N,x) \sum_{l+m+p=k}^{} (-1)^{l+m}\frac{2^{i-N-m+1}}{p!}\binom{l+i-1}{l} \right. \\
         & \ \ \ \ \ \ \left.\binom{N+m-i-1}{m}Ch_p(x)\right\}\frac{t^k}{k!}.
\end{split}\end{equation}
By comparing the coefficients on the both sides of (2.27) and (2.28), we obtain the following theorem

\begin{thm}\label{Theorem 2.2}
For $k, N=0, 1, 2, \cdot \cdot \cdot,$ we have
\begin{equation*}\begin{split}
Ch_{k+N}(x)&=k!\sum_{i=0}^{N}a_i(N,x) \sum_{l+m+p=k}^{} (-1)^{l+m}\frac{2^{i-N-m+1}}{p!}\binom{l+i-1}{l} \\
           & \ \ \ \  \times \binom{N+m-i-1}{m}Ch_p(x)\\
           &=\sum_{i=0}^{N}a_i(N,x)\sum_{l+m+p=k}^{}(-1)^{l+m}2^{i-N-m+1} \\
           & \ \ \ \  \times \binom{k}{l,m,p}(i+l-1)_l(N+m-i-1)_m Ch_p(x),
\end{split}\end{equation*}
where
\begin{equation*}\begin{split}
a_0(N,x)&=(-1)^N N!,\\
a_j(N,x)&=(x)_j(-1)^{N-j}(N-j)!\\
&\ \ \ \ \ \times\sum_{i_{j-1}=0}^{N-j}\sum_{i_{j-2}=0}^{N-j-i_{j-1}} \cdot \cdot\cdot \sum_{i_1=0}^{N-j-i_{j-1}-\cdot\cdot\cdot -i_{2}}(N-i_{j-1}-\cdot\cdot\cdot -i_{1}-j+1), \\
& \ \ \ \ \ \ \ \ \ \ \ \ \ \ \ \ \ \ \ \ \ \ \ \ \ \ \ \ \ \ \ \ \ \ \ \ \ \ \ \ \ \ \ \  \ \ \ \ \ \ \ \ \ \ \ \ \ \ \ \ \ \ \ \ \ (1 \leq j \leq N).
\end{split}\end{equation*}
\end{thm}

\end{document}